\documentstyle{article}

\begin{document}

GOLDBACH`S RULE

by Metin Aktay with thanks to Clifford H. Taubes for considering the preposterous\\

Content: 	
                               \hspace*{2.0cm}A. Overview\\
		\hspace*{2.0cm}B. List of Definitions and Variables\\
\hspace*{2.0cm}C. The Proof of Goldbach`s Conjecture\\
\hspace*{2.0cm}D. Evaluation of Observations and Further Thought \\
\hspace*{2.0cm}Appendix A \vspace*{1.0cm}

{\large \bf A. OVERVIEW }\\

Goldbach`s Conjecture, "every even number greater than $2$ can be expressed as the sum of two
 primes" is renamed Goldbach`s Rule for it can not be otherwise.

The conjecture is proven by showing that the existence of prime pairs adding to any even 
number greater than $2$ is a natural by-product of the existence of the prime sequence less
 than that even number. First it is shown that the remainder of cancellations process which 
identifies primes less than an even number also remainders prime pairs adding to that even 
number as a natural part of the process. Then a minimum limit for the number of remaindered
 prime pairs 
adding to an even number is expressed in terms of that even number  and shown to exist for 
every even number greater than $2$. Furthermore, the reasonings and formulations used in the
 proof  are demonstrated to hold against observations.\\

{\large\bf  B. LIST OF DEFINITIONS AND VARIABLES}\\
Let $ E$ be any even number $ >2$.\\
Let  $N$ be any positive integer $< E.$\\
Let $ i,n$ be counters, each of integers in natural order beginning with $1$.\\
Let a number couple be a couple where the order of two integers matters.\\
Let a number pair be a pair where the order of two integers does not matter.\\
Let $N$ symmetric integers be two integers having identical absolute difference with $N$.\\
Let $P(i)$ be any prime $\leq \sqrt{(E-1)}$ in natural order where $P(1)=2$.\\
Let $P(m)$  be the largest prime $\leq\sqrt{(E-1)}$.\\
Let $ P(i)$-prime be indivisibility by $P(i)$ where $P(i)$ is not indivisible and $1$ is indivisible.\\
Let $P(i)$-composite be divisibility by $P(i)$ where $P(i)$ is divisible.\\
Let $G1$ be the number $E/2$ symmetric $N$ which are $P(i)$-prime for all $P(i)$.\\
Let $ G2$ be the number of $E/2$ symmetric primes adding to $E.$\\
Let $GP$ be the number of prime pairs adding to $E$.\\
Let $r.f.P(i)$ be the $E/2$ symmetric $P(i)$-prime remaindering frequency for $E$ of any $P(i)$ divisor.\\
Let $stsp.m()$ be the step truncated series product with steps from $i=1$ to $i=m.$\\
Let $ PE$ be the largest prime $<E$.\\
Let $NPE$ be the number of primes $\leq PE$.\\
Let $GR$ be Goldbach Ratio,$ (GP/NPE) .$\\

{\large\bf C. THE PROOF OF GOLDBACH`S CONJECTURE}\\

{\large\bf C.1 The "Remainder" Nature Of The Prime Sequence}\\

Primeness of an integer is divisibility by no other than unity and itself. The prime sequence is identified by cancelling divisibilities and retaining indivisibilities. 
The integers not cancelled as divisible constitute the prime sequence. It is crucial to note that this identification process is an indirect process rather than direct, that
 primes are remainders, not direct creations but remnants after cancellations. \\

{\large\bf C.2 Method Of Identification Of The Prime Sequence Less Than $ E$}\\

The Sieve of Eratosthenes identifies the prime sequence up to any integer by cancelling divisibilities by prime divisors less than the square root of 
that integer and thus "remainder"ing indivisibilities. This suffices because a prime larger than that square root is multiplied by a prime less than that square 
root to produce any composite less than that integer. Therefore the sequence of primes less than $E$ are those $N$ which are not divisible by the divisors $P(i)$ 
dividing with frequency $P(i)$, where $P(i)$ were defined to be primes $\leq \sqrt{(E-1)}$.\\

{\large\bf C.3 Concepts Of $E/2$ Symmetricity And Asymmetricity of Primes and Composites}\\

On the integer line of $1$ to $(E-1)$ of $ N$, the sum of every $E/2$ symmetric integer couple is $E$.  Both members of any such couple may be prime, or both may 
be composite, or the larger member prime and the smaller composite, or the other way around. \\

A prime is defined an $E/2$ symmetric prime if it`s $E/2$ symmetric $N$ is also a prime, and an $E/2$ asymmetric prime if it`s $E/2$ symmetric $N$ is a composite. \\

A composite is defined an $E/2$ symmetric composite if it`s $E/2$ symmetric $N$ is also a composite, and an $E/2$ asymmetric composite if it`s $E/2$ symmetric $N$ is a prime.\\ 

The above may be visualised as matched integers on two integer lines matched head to tail, $(E-1)$ to $1$ and $1$ to $(E-1)$. These matched integers constitute 
couples adding to $ E$. If both members of such couples are prime then they are $ E/2$ symmetric primes, if both are composites then they are $E/2$ symmetric composites. 
Couples with one prime and one composite members contain $E/2$ asymmetric primes and $E/2$ asymmetric composites.  \\

{\large\bf C.4 Concepts of $E/2$ Symmetricity And Asymmetricity Of  $P(i)$ Divisors}\\

Any $P(i)$ divisor dividing $N$  divides either symmetrically or asymmetrically with respect to $E/2$. $E/2$ symmetric $P(i)$ divisors divide 
$E/2$ and $(E/2)+nP(i)$ and $(E/2)-nP(i)$. $E/2$ asymmetric $P(i)$ divisors never divide $E/2$ symmetric $N$.\\

$E/2$ symmetricity of any $P(i)$ divisor depends on divisibility of $E/2$ by that $P(i)$. If $E/2$ is divisible by a $P(i)$ then that $P(i)$ divides and  "remainder"s 
symmetrically with respect to $E/2$. If $E/2$ is not divisible by a $P(i)$ then that $P(i)$ divides asymmetrically with respect to $E/2$. This dependence applies to 
all $P(i)$ except $P(1)=2$ which is $E/2$ symmetric independent of $E/2$ divisibility, with either two divisibilities or two indivisibilities bracketing $E/2$.\\

If all $P(i)$ divisors were $E/2$ symmetric, then both composite and prime $N$ would be $E/2$ symmetric, and thus all primes $<E$ would be members of 
 prime couples adding to $E$, except for $P(i)$ themselves which are symmetric with composites divided by themselves. If all $P(i)$ divisors were to be $E/2$ asymmetric, 
except for $P(1)=2$ which can not be so, then there would be minimal $E/2$ symmetric composites and primes.\\

{\large\bf  C.5 Concept Of $P(i)$-Primes And The Primes As Intersecting Sets Of $P(i)$-primes}\\

If being $P(i)$-prime is indivisibility by $P(i)$, where $1$ is defined $P(i)$-prime and $P(i)$ itself is defined divisible, there would be $1$ 
$P(i)$-composite and $(P(i)-1)$  $ P(i)$-primes every $P(i)$ consecutive $N$ for any $P(i)$ divisor. The $N$ which are $P(i)$-prime for all $P(i)$ can not be but prime 
except for $1$. Therefore, the primes $<E$ are the $N$ which are $P(i)$-prime for all $P(i)$, with the addition of the  $P(i)$ themselves
the deduction of $1$.\\

\newpage
{\large\bf C.6 Existence And Frequency Of $E/2$ Symmetric $P(i)$-primes For Any $P(i)$}\\

Any $E/2$ symmetric $P(i)$ divisor will divide $1$  $E/2$ symmetric $P(i)$-composite and remainder $(P(i)-1)$  $ E/2$ symmetric $P(i)$-primes every $P(i)$ consecutive $N.$\\

Any $E/2$ asymmetric $P(i)$ divisor will divide $1$  $E/2$ asymmetric $P(i)$-composite, and will remainder $1$ $E/2$ 
asymmetric $P(i)$-prime and  $(P(i)-2)$  $ E/2$ symmetric $P(i)$-primes
 for every $P(i)$ consecutive $N$.\\ 

The above may be visualised with the two lines of integers matched head-to-tail. Any $P(i)$ divisor cancelling $P(i)$-composites 
and remaindering $P(i)$-primes on each line begins dividing from opposite ends, thus may or may not meet at the midpoint since that $P(i)$ may
 or may not be $E/2$ symmetric. Let the two lines be counted by $n$ counting from one end only, beginning with $n=1$, which counts one line forward and 
the other line backwards. If the cancellations by any $P(i)$ divisor on each line is accounted for with respect to $n$, then the frequency with respect to $n$ of divisibility 
by that $P(i)$ of  both lines is $1/P(i)$ since both lines move by $1$ for every unit change in $n$. However, for an $E/2$ asymmetric $P(i)$ divisor, the divisibility of each line 
by that $P(i)$ may lead or lag the other line with respect to $n$. If a $P(i)$ divisor is $E/2$ symmetric it will divide one matched integer couple and remainder $((P(i)-1)$ matched 
integer couples as $P(i)$-primes every $P(i)$  $n$. If a $P(i)$ divisor is $E/2$ asymmetric it will divide $1$  $P(i)$-composite on the first line matched with a $P(i)$-prime on the second, and 
$1$ $P(i)$-composite on the second line matched with a  $P(i)$-prime on the first, and will remainder $(P(i)-2)$ couples with $P(i)$-primes on each line every $P(i)$  $n$.\\

$P(1)=2$ divisor can not be $E/2$ asymmetric. $P(1)=2$ would render all $2$-Primes $E/2$ asymmetric if it could be $E/2$ asymmetric.  $P(1)=2$ divisor will always 
remainder $1$ $E/2$ symmetric  $2$-Prime for every $ 2$ consecutive $N$.\\

$E/2$ asymmetric $P(2)=3$ divisor will remainder $1$  $E/2$ symmetric $3$-prime for every $3$ consecutive $N$ and $E/2$ symmetric $P(2)=3$ will remainder $2$. $P(i)$ divisors 
larger than $3$ will remainder more than $1$  $E/2$ symmetric $P(i)$-primes for every  $P(i)$ consecutive $N$, i.e. $E/2$ asymmetric $ P(3)=5$ 
divisor will remainder $3$  $E/2$ symmetric $5$-primes 
every $5$ consecutive $N$ and $E/2$ symmetric $P(3)=5$ will remainder $4$.\\
 
Therefore $E/2$ symmetric $P(i)$-primes exist for each $P(i)$ divisor of  $ E$.\\
\newpage
{\large \bf C.7 Lower Limit For The Number of $E/2$ symmetric $P(i)$-Primes for all $P(i)$}\\

To find the number of $E/2$ symmetric $P(i)$-primes for a $P(i)$ divisor, the number of  $N$, which is $E-1$, 
is multiplied with the $E/2$ symmetric
 $P(i)$-prime remaindering frequency of that $P(i)$. Therefore \\
\hspace*{2.0cm}$number \ \ of\ \  E/2\ \  symmetric \ \ P(i)-primes = (E-1)*r.f.P(i)$\\
This may result in a fractional result which may need truncation or rounding up. Therefore\\
\hspace*{2.0cm}$number\ \ of\ \ E/2\ \  symmetric \ \ P(i)-primes \geq  truncate ((E-1)*r.f.P(i))$\\

To find a lower limit for $G1$, the number of $N$ which are $E/2$ symmetric $P(i)$-prime for all $P(i)$, $E-1$ 
is multiplied consecutively by $r.f.P(i)$ with truncations at each step. This reduction with consecutive 
remaindering frequencies is valid because $P(i)$ are indivisible by each other.\\
Therefore \\
\hspace*{2.0cm}$G1\geq  stsp.m ( (E-1)* r.f.P(i) )$\\
Expectation of $G1$ will be minimum where all $r.f.P(i)$ are minimum.\\ 
Therefore\\
\hspace*{2.0cm}$min.G1 \geq  stsp.m ( (E-1)*sp(min.r.f.P(i)))$.\\
All $r.f.P(i)$ will be minimum where all $P(i)$ divisors are $E/2$ asymmetric. \\
Therefore\\
\hspace*{2.0cm}$min.r.f.P(i) = ((P(i)-2)/P(i))$, except for $r.f.P(1)=1/2$.\\

If any product $( (a)* (b/c) )$, where $a,b,c$ are positive integers, and where $c \leq a$, is truncated, then 
the result can not be less than $b$. With this logic, lower limits can be found for the successive truncations
 of the step truncated series product by comparing the denominator of each multiplier with the numerator of 
the multiplied.\\

Let $stsp.m ( (E-1)*sp(min.r.f.P(i)))$ begin from $min.r.f.P(m)$ and work backwards, and let the first
 truncation result be $T(1)$ and the last truncation result be $T(m)$, where  $T(m) \leq  min.G1.$\\ 
Therefore\\
\hspace*{2.0cm}$T(1) = truncate ( (E-1) * min.r.f.(P(m)) )$.\\
Given that except for $P(1)=2$\\
\hspace*{2.0cm}$min.r.f.P(i)= ( (P(i)-2)/P(i) )$\\ 
Then\\
\hspace*{2.0cm}$T(1) = truncate ( (E-1) * ((P(m)-2)/P(m)) )$\\
Given that\\
\hspace*{2.0cm}$	(E-1) \geq  (P(m))^2$\\
Then\\
\hspace*{2.0cm}$	T(1) \geq  ( P(m) * (P(m)-2) ).$\\
Given that \\
\hspace*{2.0cm}$	T(2) = truncate ( T(1) * min.r.f.(P(m-1)) )$\\
and \\
\hspace*{2.0cm}$	 min.r.f.P(m-1)= ( (P(m-1)-2)/P(m-1) ) $\\
Then\\
\hspace*{2.0cm}$T(2) = truncate ( T(1) * (P(m-1)-2) / P(m-1) )$\\
Then\\
\hspace*{2.0cm}$T(2) \geq truncate ( P(m) * (P(m)-2) * (P(m-1)-2) / P(m-1) )$\\ 
Given that \\
\hspace*{2.0cm}$	P(i) \leq  (P(i+1)-2)$ except for $P(1)=2$\\
Then\\
\hspace*{2.0cm}$ 	P(m-1)\leq  (P(m)-2) $\\
Therefore\\
\hspace*{2.0cm}$	T(2)\geq  ( P(m) * (P(m-1)-2) )$ .\\
If lower limits for $T(n)$ are found consecutively as above\\ 
Then \\
\hspace*{2.0cm}$	T(n)\geq ( P(m) * ( P(m-n+1)-2))$.\\
Therefore\\ 
\hspace*{2.0cm}$	T(m-1)\geq  ( (P(m)) * (P(2)-2) )$.\\  
Since $P(2)=3$ then\\
\hspace*{2.0cm}$	T(m-1) \geq  ( (P(m) )$\\
Given that $r.f.P(1)=1/2$, then\\
\hspace*{2.0cm}$	T(m) \geq  truncate  ( (P(m)) * 1/2 )$\\
Therefore\\
\hspace*{2.0cm}$min.G1\geq  truncate ( P(m) /2 )$.\\

{\large \bf C.8 Lower Limit For The Number Of Prime Pairs Adding To $E$}\\

Since $1$ is defined $P(i)$-prime for all $P(i)$, $G1$, which is the number of $N$  $E/2$ symmetric 
$P(i)$-prime for all $P(i)$, may count $1$ if the $E/2$ symmetric counterpart of $1$ is also $P(i)$ prime for 
all $P(i)$. Given that $G2$ is the number of $E/2$ symmetric primes adding to $E$, and given that $G1$ may count 
$1$ and $(E-1)$, and given that $G1$ excludes $P(i)$, and given that $min.G1$ is a truncated result, then\\
$G2 \geq  min.G1 -( 2\  or\  0) +  2*(number\ \  of\ \  P(i)\ \  adding \ \ to\ \ E \ \ with\ \ 
 other\ \ primes) $\\
Therefore\\
\hspace*{2.0cm}$	min.G2 \geq (min.G1 -2)$.\\
Therefore\\
\hspace*{2.0cm}$	min.G2\geq (truncate(P(m)/2) -2)$\\ 
$GP$, the number of prime pairs adding to $E$, is half of $G2$ since every $E/2$ symmetric prime is a member 
of a pair of primes adding to $E$. If the halving of $G2$ is fractional, then it is rounded up since an odd 
numbered $G2$ indicates that $E/2$ itself is prime and is counted once. $G2$ is even numbered where $E/2$ 
itself is not prime. \\
 
Therefore\\
\hspace*{2.0cm}$	min.GP\geq  rounded-up (1/2*min.G2)$.\\
Then\\
\hspace*{2.0cm}$	min.GP\geq  rounded-up (1/2* ( truncate(P(m)/2) -2 ) ) $.\\
Therefore\\
\hspace*{2.0cm}$	min.GP \geq 1$ for $E\geq 50$ where $P(m)\geq 7$.\\

The logic of the above proof also proves that there is at least one prime pair adding to 
$E$ with primes $>P(m)$ for every $E \geq  50$.\\
 
Given the existence of at least one prime pair adding to $E$ for every $E<50$, then\\
\hspace*{2.0cm}$	min.GP\geq  1$ for every $E$.\\
Therefore\\
\hspace*{2.0cm}$	GP\geq  1$ for every $E$.\\
Magnitude of $min.GP \cong \sqrt{E}/4$ for large $E$ since magnitude of $P(m)\cong \sqrt{E}$ for large $E$.\\

{\large\bf C.9  Conclusion Of The Proof of Goldbach`s Conjecture}\\ 

Thus it is proven that the remaindering of the prime sequence $<E$ by divisors $P(i)\leq\sqrt{(E-1)}$ can 
not avoid remaindering $E/2$ symmetric primes adding to $E$ for every $E$ even in minimal expectation 
conditions, since the number of prime pairs adding to $E$ is $\geq  1$ for every $E$, and since a hypothetical
 minimum limit for this number is $\cong \sqrt{E}/4$ for large $E$.\\

Therefore there will always be $E/2$ symmetric prime pairs adding to $E$ for any $E$ 
since the remaindering process for the prime sequence $<E$ can not avoid remaindering $E/2$ 
symmetric primes. Therefore there will be at least one pair of $E/2$ symmetric primes adding to 
$E$ for any $E$. \\

Therefore every even number greater than $2$ can be expressed as the sum of two primes.\\

\newpage
{\large\bf D. EVALUATION OF OBSERVATIONS AND FURTHER THOUGHT}\\

{\large\bf D.1 Table Of Observations And Calculations}\\

Behaviour with respect to $E$ of the number of prime pairs adding to $E$ may be deduced from 
Appendix A, Observations and Calculations, where the formulations utilised above and actual counts 
are tabulated for sample $E$ evenly and saliently spread up to $10,000$.\\
 
{\large\bf D.2 Validation Of Assumed Relationships}\\

$GP$ may be approximated as $(E/2)*sp(r.f.P(i))$, which would undercount since it fails to check 
the $P(i)$  excluded by $r.f.P(i)$, and would undercount since it omits rounding up where needed, 
and would overcount since it omits truncation where needed, and would overcount since it may count $1$. 
In Appendix A, $GP \cong  (E/2)*sp(r.f.P(i))$ is calculated as is without the reduction rationalisations used in
 the proof. It is observed that the proportion of error in this raw term decreases with larger $E$. More important
 than decreasing error is that it is observed to track flawlessly the volatility with respect to $E$ of actual $GP$. \\

{\large\bf D.3 Observations and Evaluations Of Goldbach Ratio}\\

Appendix A shows that $GR$ varies as expected. $GR$ is low where $E/2$ is indivisible by smaller $P(i)$
 and high where $E/2$ is divisible by smaller $P(i)$, since divisibility of $E/2$ by smaller $P(i)$ increases 
$E/2$ symmetricity of primes substantially.\\

$E$ with $E/2$ prime or divisible only by $P(1)=2$ have low $GR$. All $E$ with $E/2$ divisible by $3$ have high $GR$
 because
 avoidance of an $E/2$ asymmetric $3$ divisor doubles symmetricity of primes. For example, $E=210$ has a high $GR$ , 
$41\%$, 
 since $E/2=105$ is divisible by $3, 5, 7$, which thus  remainder primes symmetrically with respect to $E/2$. A $GR$ value 
of $41\%$ is considered high since $GR$ is defined by utilising pairs of primes adding to $E$, which means that the maximum 
possible $GR$ is $50\%$ , where every prime is a member of a pair of primes adding to $E$.\\

{\large\bf D.4 Deductions In Relation To The Goldbach Comet}\\

$E/2$ symmetricity of $P(i)$ divisors, which is divisibility of $E/2$ by $P(i)$, explains the dense cluster 
bands which form when $GP$ is plotted against $E$, a plot called the Goldbach Comet on account of these cluster bands. 
The densest asymptotic cluster band is formed by $E$ with $E/2$ asymmetric smaller prime divisors. $E/2$ symmetricity of 
smaller prime divisors explain other bands. The next dense band are the $E$ with $E/2$ symmetric $P(2)=3$ divisor, at 
$(2/3)/(1/3)=2$ times the heights of the asymptotic lowest band. Then the next dense band are the $E$ with $E/2$ symmetric
 $P(3)=5$ divisor, at $(4/5)/(3/5)=4/3$ times the heights of the asymptotic lowest band. Joint $E/2$ symmetricity of the 
smaller primes also form distinct dense bands, such as the $E$ with $E/2$ symmetric $P(2)=3$ and $P(3)=5$ divisors, at
 $(2/3*4/5)/(1/3*3/5)=8/3$ times the heights of the asymptotic lowest band.\\

{\large\bf D.5 Further Thought}\\

The number of prime pairs adding to $E$ increases with $E$ but with high volatility. It could be that 
this volatility is fluctuation around a fundamental relationship. Assuming  the prime sequence to be a 
discrete wave function, and assuming the estimate $N/lnN$ of J.S.Hadamard to be it`s frequency for $N$, this fundamental
 relationship is likely to be $GP = E/(2*(lnE)^2)$.

\vspace{\fill}
Metin Aktay\\
Ihsan Aksoy sok. EVA apt. No:7/2,  Camlik, Etiler, Besiktas, Istanbul 80630, Turkey\\
Phone: +90 212 2651016 	Fax: +90 212 2577374 	Mobile: +90 532 2741771\\
E-mail: maktay@superonline.com, or maktay@mba1979.hbs.edu
\newpage
{\bf APPENDIX A: OBSERVATIONS AND CALCULATIONS}

$
\begin{tabular}[h]{|r|r|r|r|r|r|r|r|r|r|}   \hline
E&PE &NPE &P(m) &E/2 & E/2\ \ Factors\ \ &observed  & GR\%& calcul.  &Error \% \ \ of  \\ 
& & & &&$ \leq P(m)$ &GP&&GP&calcul.\ \ GP\\ \hline  
128&127&31&11&64&2&3&10&4&33 \\ \hline
210 &199&46&13&105&2,3,5,7&19&41&17&-11\\ \hline
222&211&47&13&111&P&11&23&5&-55\\ \hline
502&499&95&19&251&P&15&16&10&-33\\ \hline
512&509&97&19&256&2&11&11&10&-9\\ \hline
678&677&123&23&339&3&28&23&24&-14\\ \hline
1,006&997&168&31&503&P&18&11&16&-11\\ \hline
1,024&1,021&172&31&512&2&22&13&16&-27\\ \hline
1,510&1,499&239&37&755&5&33&14&30&-9\\ \hline
2,018&2,017&306&43&1,009&P&28&9&27&-4\\ \hline 
2,048&2,039&309&43&1,024&2&25&8&27&8\\ \hline
2,490&2,477&367&47&1,245&3,5&94&26&85&-10\\ \hline
3,022&3,019&433&53&1,511&P&42&10&37&-12\\ \hline
3,514&3,511&490&59&1,757&7&51&10&50&-2\\ \hline
4,006&4,003&552&61&2,003&P&52&9&46&-12\\ \hline
4,096&4,093&564&61&2,048&2&53&9&47&-11\\ \hline
4,690&4,679&633&67&2,345&5,7&95&15&83&-13\\ \hline
5,006&5,003&670&67&2,503&P&63&9&56&-11\\ \hline
5,610&5,591&738&73&2,805&2,3,5,11,17&198&27&186&-6\\ \hline
6,002&5,987&783&73&3,001&P&62&8&63&2\\ \hline
6,578&6,577&851&79&3,289&2,11,13,23&89&10&86&-3\\ \hline
7,022&7,019&903&83&3,511&P&72&8&70&-3\\ \hline
7,314&7,309&932&83&3,657&2,3,23,53&172&18&156&-9\\ \hline
8,002&7,993&1,007&89&4,001&P&80&8&78&-3\\ \hline
8,192&8,191&1,028&89&4,096&2&76&7&80&5\\ \hline
8,610&8,609&1,072&89&4,305&2,3,5,7,41&282&26&276&-2\\ \hline
9,014&9,013&1,021&89&4,507&P&96&9&88&-8\\ \hline
9,510&9,497&1,177&97&4,755&3,5&253&21&243&-4\\ \hline
9,998&9,973&1,229&97&4,999&P&99&8&96&-3\\ \hline

\end{tabular}
$
 $\ $\\
Sample calculations for $GP$ calculated as $((E-1)/2)*sp(r.f.P(i))$\\
$GP$ for $E=2490: \ (2,489/2)*(1/2*2/3*4/5*5/7*9/11*11/13*15/17*17/19*21/23*27/29*29/31*35/37*39/41*41/43*45/47)$\\
$GP$ for $E=3022: \ (3021/2)*(1/2*1/3*3/5*5/7*9/11*11/13*15/17*17/19*21/23*27/29*29/31*35/37*39/41*41/43*45/47*51/
53)$

\end{document}